\documentclass[11pt]{article}
\usepackage{amsmath,amsfonts,latexsym,graphicx,amssymb,amsthm}
\usepackage[tips,matrix,arrow]{xy}
\usepackage{etex}
\usepackage{tikz}

\usepackage{xcolor}

\definecolor{violet}{cmyk}{1,0,0,0}

\usepackage{hyperref}

\usepackage{setspace}

\newcommand{\R}{\mathbb{R}}
\newcommand{\Z}{\mathbb{Z}}

\newcommand{\Aut}{\mathrm{Aut}}

\newtheorem{thm}{Theorem}[section]
\newtheorem{defn}{Definition}[section]

\newtheorem{ack}{Acknowledgment}

\newcommand{\vep}{\varepsilon}
\newcommand{\rad}{\mathrm{rad}}

\begin{document}
\title{Classification of marked elliptic root systems \\ with non-reduced affine quotient}
\author{A. Fialowski, K. Iohara and Y. Saito}

\maketitle
\begin{abstract}
The class of root systems, called elliptic root systems, were introduced in 1985 by K. Saito, for his studies on a normal surface singularity which contains a regular elliptic curve in its minimal resolution. He also classified such root systems when they admit a reduced affine quotient, as root system. In this note, we provide the classification of elliptic root systems that admit a non-reduced affine quotient, thus complete the classification of such root systems. 
\end{abstract}

\begin{center}\begin{flushleft}
{\small{{\bf Résumé}~ La classe des systèmes de racines, dits elliptiques, a été introduite par K. Saito en 1985 pour ses études sur une singularité normale d'une surface dont sa résolution minimale contient une courbe elliptique régulière. Il a même classifié tels systèmes de racines dans le cas où ils admettent un quotient affine réduit, comme système de racines. Dans cette note, nous allons donner la classification de systèmes elliptiques de racines qui admettent un quotient affine non-réduit, d'où nous achevons la classification des systèmes elliptiques de racines.}}
\end{flushleft}\end{center}
\begin{center}
\end{center}

\pagestyle{plain}
\section{Introduction}
In his study on simply elliptic singularities (cf. \cite{Saito1974}), K. Saito introduced the notion of \textbf{elliptic root systems}, what were called \textit{$2$-extended affine root systems}, in 1985 \cite{Saito1985}. 
Such root systems are defined in a real vector space $F=\R^{l+2}$ ($l>0$) equipped with a positive semi-definite metric $I$ whose radical $\rad(I)$ is of dimension $2$. Hence, it can be viewed as a $2$-dimensional generalization of affine root systems. However, this newly defined root system has an interesting property: it has a finite order Coxeter element in its Weyl group !  This fact plays a crucial role for its study and has important geometric applications in singularity theory.  \\

Let $R$ be an elliptic root system and let $G$ be a one-dimensional subspace of the radical $\rad(I)$ of $I$ such that the sublattice $G \cap Q(R)$ of $Q(R):=\Z R$ is full in $G$. K. Saito has classified the pairs $(R,G)$ of an elliptic root system $R$ with its marking $G$, under the assumption that the quotient $R/G$, which is the image of $R$ in $F/G$ via the canonical projection $F \twoheadrightarrow F/G$, is a reduced affine root system. Hence, his classification heavily depends on the classification of affine root systems due to V. G. Kac \cite{Kac1969}, R. V. Moody \cite{Moody1969} and I. G. Macdonald \cite{Macdonald1972}, where Kac and Moody considered affine root systems as the root system of affine Lie algebras. Notice that Macdonald classified affine root systems without discussing their relations with Lie algebras, hence including non-reduced affine root systems. 
Saito's classification doesn't imply the classification of reduced marked elliptic root systems $(R,G)$, since $R/G$ can be non-reduced even if $R$ itself is reduced. \\

In this note, we give the classification theorem of the marked elliptic root systems $(R,G)$ whose quotient $R/G$ is a non-reduced affine root system, thus complete the classification of marked elliptic root systems. 

\medskip
\begin{ack}{\rm
We want to thank the referee for pointing out the paper \cite{AKY} which, together with \cite{AABGP}, helped us find a missing case in the classification.
Research of Y. S.  is supported by JSPS KAKENHI Grant number JP20K03568.
}\end{ack}

\section{Marked elliptic root systems}
Let $F$ be a real vector space and $I: F \times F \rightarrow \R$ be a symmetric bilinear form whose signature is $(l_+, l_0, l_-)$ for some non negative integers $l_+, l_0, l_- $ such that not all of them are zero. As usual, for any non-isotropic vector $\alpha \in F$, we set
\[ \alpha^\vee=\frac{2}{I(\alpha,\alpha)}\alpha \in F, \]
and define an isometry $w_\alpha\in O(F,I)$ by
\[ w_\alpha(\lambda)=\lambda-I(\lambda, \alpha^\vee)\alpha. \]
\begin{defn} A non-empty discrete subset $R$ of $F$ is called a \textbf{generalized root system} if it satisfies
\begin{enumerate}
\item The lattice $Q(R)$ (called the \textbf{root lattice}), spanned by the elements of $R$, is full in $F$, i.e., $\R \otimes_{\Z} Q(R) \cong F$. 
\item For any $\alpha \in R$, one has $I(\alpha, \alpha) \neq 0$. 
\item For any $\alpha, \beta \in R$, one has $I(\alpha^\vee, \beta) \in \Z$.
\item For any $\alpha \in R$, one has $w_\alpha(R)=R$.
\item Assume that there exists two subsets $R_1, R_2$ of $R$ which are orthogonal and $R_1 \cup R_2=R$, then either $R_1$ or $R_2$ is empty. 
\end{enumerate}
\end{defn}
\vskip 0.1in

An \textbf{elliptic root system} $R$ is a generalized root system belonging to $F$ with a metric $I$ whose signature is $(l,2,0)$. A vector subspace $G$ of $\rad(I)$ is called a \textbf{marking} if the sublattice $G \cap Q(R)$ of $G$ is full in $G$. The pair $(R,G)$ of an elliptic root system with its marking is called a \textbf{marked elliptic root system}, \textit{mERS} for short. The image of $R$ by the canonical projection $F \twoheadrightarrow F/G$ is known to be an affine root system, denoted by $R/G$ and is called the \textbf{quotient root system.}
\\

By definition, the sublattice $\rad_{\Z}(I):=Q(R) \cap \rad(I)$ of $Q(R)$
is of rank $2$. Hence, there exists two non-zero vectors $a$ and $b$ of $F$ which generate the lattice $\rad_{\Z}(I)=\Z a \oplus \Z b$. Here and after, we always fix $G=\R a$. 
In 1985, K. Saito \cite{Saito1985} classified mERSs $(R,G)$ whose quotient $R/G$ is a reduced affine root system. Here, a root system is said to be \textbf{reduced} if for any $\alpha \in R$, neither $2\alpha$ nor $\frac{1}{2}\alpha$ belongs to $R$. 

\section{Main results}
In the rest of this note, we fix an ambient real vector space $F=\R^{l+2}$ equipped with a metric $I$ whose signature is $(l,2,0)$. For each root system $R$, we denote the set of short (resp. middle length, long) roots by $R_s$, $R_m$ and $R_l$, respectively. 

\subsection{New root systems}
First, we introduce some reduced mERSs $(R,G)$ with non-reduced affine quotient $R/G$. Since there are only 6 such cases, we call each type of the marked root system by making reference to $BC_l$. 
\subsubsection{Reduced classical type}
\label{sect:intro_red-classic}
Each mERS $(R,G)$ with $G=\R a$ is defined by
\[ R=(R(X_l)_s+\Z a) \cup (R(X_l)_m+\Z a) \cup (R(X_l)_l+(1+2\Z)a), \]
where $R(X_l)$ is a non-reduced affine root system, and for each $R$, it is given by the following table: \\
\begin{center}
\begin{tabular}{|c||c|c|c|c|} \hline
Type of $R$ & $BC^{(1,2)}_l$ & $BC_l^{(4,2)}$ & $BC_l^{(2,2)\sigma}(1)$ & $BC_l^{(2,2)\sigma}(2)$ 
 \\ \hline
$X_l$ & $BCC_l$ & $C^\vee BC_l $ & $BB_l^\vee$ & $C^\vee C_l$  \\ \hline
\end{tabular}
\end{center}
Here we recall that these $4$ non-reduced affine root systems are defined as follows (cf. \cite{Macdonald1972}):
\begin{align*}
R(BCC_l)=
&R(BC_l)+\Z b 
\phantom{ \;\, R \cup(R(BC_l)_m+2\Z b) \cup (R(BC_l)_l+{\color{red}4}\Z b)\qquad} 
(l\geq 1), \\
R(BB_l^\vee)=
&(R(BC_l)_s+\Z b) \cup (R(BC_l)_m+\Z b) \cup (R(BC_l)_l+2\Z b) \qquad 
\;\, \,(l\geq 2), \\
R(C^\vee C_l)=
&(R(BC_l)_s+\Z b) \cup (R(BC_l)_m+2\Z b) \cup (R(BC_l)_l+2\Z b) \qquad
(l\geq 1), \\
R(C^\vee BC_l)=
&(R(BC_l)_s+\Z b) \cup (R(BC_l)_m+2\Z b) \cup (R(BC_l)_l+4\Z b) \qquad
(l\geq 1),
\end{align*}
and the root system of type $BC_l$ is defined by 
\[ R(BC_l)=\{ \, \pm \vep_i\, \}_{1\leq i\leq l} \cup \{ \, \pm (\vep_i\pm \vep_j)\, \vert \, 1\leq i<j\leq l\, \}
\cup \{ \, \pm 2\vep_i\, \}_{1\leq i\leq l}, \]
where $\{\vep_i\}_{1\leq i \leq l}$ satisfy $I(\vep_i,\vep_j)=\delta_{i,j}$, the Kronecker delta. 
\subsubsection{Reduced $\ast$-type}\label{section:intro_red-star}
There are two ERSs of $\ast$-type, and they are defined by
\begin{align*}
R(BC_l^{(1,1)\ast})=
&(R(BCC_l)_s+\Z a) \cup (R(BCC_l)_m+\Z a) \cup (R(BC_l)_l+L_{1,1}), \\
R(BC_l^{(4,4)\ast})=
&(R(BC_l)_s+L_{1,1}) \cup (R(C^\vee BC_l)_m+2\Z a) \cup (R(C^\vee BC_l)_l+4\Z a),
\end{align*}
where we set 
\[ L_{1,1}=\{ \, ma+nb\, \vert \, (m-1)(n -1) \equiv 0 [2]\, \}. \]
Notice that they had been discovered by S. Azam \cite{Azam2002} in 2002.  \\

Second, we consider the non-reduced root systems. Since there are many such root systems, we call each type of the marked root systems by making reference to its affine quotient $R(X_l)=R/G$ which is one of the type $BCC_l, C^\vee BC_l, BB_l^\vee$ and $C^\vee C_l$. 
\subsubsection{Non-reduced classical type}\label{sect:intro_non-red-classic}
We introduce $4$ types of the root systems: $X_l^{(1)}, \, X_l^{(2)}(1), \,X_l^{(2)}(2)$ and 
$X_l^{(4)}$. They are defined by
\begin{align*}
R(X_l^{(1)})=
&R(X_l)+\Z a, \\
R(X_l^{(2)}(i))=
&(R(X_l)_s+\Z a) \cup (R(X_l)_m+i \Z a) \cup (R(X_l)_l+2\Z a) \qquad i=1,2, \\
R(X_l^{(4)})=
&(R(X_l)_s+\Z a) \cup (R(X_l)_m+2\Z a) \cup (R(X_l)_l+4\Z a). 
\end{align*}
Here, $R(X_l)$ is one of the $4$ types: $BCC_l, C^\vee BC_l, BB_l^\vee$ and $C^\vee C_l$.
\subsubsection{Non-reduced $\ast$-type}\label{sect:intro_non-red-star}
In this case, there can be several $\ast$-types : $\ast_i, \, \ast_{i'}$ ($i=0,1$), $\ast_s$ and $\ast_l$. 

Let $L_{i,j}, L_{i,j}^{s_1,s_2}$ ($i,j=0,1$) and ($s_1, s_2 \in \Z_{>0}$) be the subsets of the lattice $\rad_{\Z}(I)=\Z a\oplus \Z b$ defined as follows ($L_{1,1}$ has already been introduced):
\begin{align*}
L_{i,j}=
&\{\, ma+nb\, \vert \, (m-i)(n-j)\equiv 0 [2]\, \}, \\
L_{i,j}^{s_1,s_2}=
&\{\, s_2ma+s_1nb\, \vert \, (m-i)(n-j)\equiv 0 [2]\, \}. 
\end{align*}
Let $t_1$ be an integer defined by the following table:
\begin{center}
\begin{tabular}{|c||c|c|c|} \hline
$X_l$ & $BCC_l$ & $C^\vee BC_l $ & $C^\vee C_l$  \\ \hline
$t_1$ & $1$ & $4$ & $2$  \\ \hline
\end{tabular}
\end{center}
and $t_2$ an integer such that $(t_1, t_2) \in \{1,2,4\}^2 \setminus \{ (1,4), (4,1)\}$.  \\

The set $R_m$ always is of the form
\[ R_m=R(X_l)_m+\min\{2,t_2\}\Z a. \]
and the set $R_s$ and $R_l$ are described as follows: \\
\vskip 0.2in
\noindent{\fbox{\textbf{$\ast_i$-type ($i=0,1$)}}} \hskip 0.1in 
For each such pair $(t_1,t_2)$, the subsets $R_s$ and $R_l$  are given by
\begin{enumerate}
\item for $i=0$ and
\begin{enumerate}
\item[i)] $t_1, t_2 \in \{1,2\}$ and $(t_1, t_2) \neq (2,2)$, 
\[ R_s=R(X_l)_s+\Z a \qquad  \text{and} \qquad R_l=R(BC_l)_l+L_{0,0}^{t_1,t_2}, \]
\item[ii)] $t_1=t_2=2$, 
\[ R_s=R(BC_l)_s+L_{0,0} \qquad \text{and} \qquad R_l=R(BC_l)_l+L_{0,0}^{2,2}, \]
\item[ii)] $t_1,t_2 \in \{2,4\}$ and $(t_1, t_2) \neq (2,2)$, 
\[ R_s=R(BC_l)_s+L_{0,0} \qquad \text{and} \qquad R_l=R(X_l)_l+t_2\Z a,\]
\end{enumerate}
\item for $i=1$ and 
\begin{enumerate}
\item[i)] $t_1, t_2 \in \{1,2\}$ and $(t_1, t_2) \neq (2,2)$, 
\[ R_s=R(X_l)_s+\Z a \qquad \text{and} \qquad R_l=R(BC_l)_l+L_{1,1}^{t_1,t_2}, \]
\item[ii)] $t_1=t_2=2$,
\[ R_s=R(BC_l)_s+L_{0,0} \qquad \text{and} \qquad R_l=R(BC_l)_l+L_{1,1}^{2,2}, \]
\item[iii)] $t_1,t_2 \in \{2,4\}$ and $(t_1, t_2) \neq (2,2)$, 
\[ R_s=R(BC_l)_s+L_{1,1} \qquad \text{and} \qquad R_l=R(X_l)_l+t_2\Z a. \]
\end{enumerate}
\end{enumerate}
For $i=0,1$, the nomenclature of each mERS given by the above formula is as follows:
{\small
\begin{center}
\begin{tabular}{|c||c|c|c|c|c|c|c|} \hline
$(t_1,t_2)$ & $(1,1)$ & $(1,2)$ & $(2,1)$ & $(2,2)$ & $(2,4)$ & $(4,2)$ & $(4,4)$\\ \hline
Type of $(R,G)$ & $BCC_l^{(1)\ast_i}$ & $BCC_l^{(2)\ast_i}$ & $C^\vee C_l^{(1)\ast_i}$  
& $C^\vee C_l^{(2)\ast_i}$ & $C^\vee C_l^{(4)\ast_i}$ & $C^\vee BC_l^{(2)\ast_i}$ & $C^\vee BC_l^{(4)\ast_i}$\\ \hline
\end{tabular}
\end{center}}
Remark that the ERSs
$R(BCC_l^{(1)\ast_1})$ and $R(C^\vee BC_l^{(4)\ast_1})$ are reduced and they are exactly the ERSs $R(BC_l^{(1,1)\ast})$ and $R(BC_l^{(4,4)\ast})$, respectively. 
\vskip 0.2in

\noindent{\fbox{\textbf{$\ast_{i'}$-type ($i=0,1$)}}} \hskip 0.1in 
One has $t=t_1=t_2 \in \{1,2,4\}$, and the subsets $R_s$ and $R_l$  are given by
\begin{enumerate}
\item for $t=1$,
\[ R(BCC_l^{(1)\ast_{0'}})_s=R(BCC_l)_s+\Z a \qquad \text{and} \qquad R(BCC_l^{(1)\ast_{0'}})_l=R(BC_l)_l+L_{0,1}, \]
\item for $t=2$,
\[ R(C^\vee C_l^{(2)\ast_{1'}})_s=R(BC_l)_s+L_{0,0} \qquad \text{and} \qquad 
R(C^\vee C_l^{(2)\ast_{1'}})_l=R(BC_l)_l+L_{1,0}^{2,2} \]
\item for $t=4$, 
\[ R(C^\vee BC_l^{(4)\ast_{0'}})_s=R(BC_l)_s+L_{0,1} \qquad \text{and} \qquad
   R(C^\vee BC_l^{(4)\ast_{0'}})_l=R(C^\vee BC_l)_l+4\Z a.  \]
\end{enumerate}
\vskip 0.2in 

\noindent{\fbox{\textbf{$\ast_{\natural}$-type ($\natural \in \{s,l\}$)}}} \hskip 0.1in 
One has $t_1=t_2=2$, and the subsets $R_s$ and $R_l$  are given by
\begin{enumerate}
\item for $\natural=s$,
\[ R(C^\vee C_l^{(2)\ast_{s}})_s=R(BC_l)_s+L_{0,0}  \qquad \text{and} \qquad
   R(C^\vee C_l^{(2)\ast_{s}})_l=R(C^\vee C_l)_l+2\Z a, \]
\item for $\natural=l$,
\[  R(C^\vee C_l^{(2)\ast_{l}})_s=R(C^\vee C_l)_s+\Z a \qquad \text{and} \qquad 
    R(C^\vee C_l^{(2)\ast_{l}})_l=R(BC_l)_l+L_{0,0}^{2,2}. \]
\end{enumerate}
\vskip 0.2in
Remark that the numbers $t_i$ ($i=1,2$) are the so-called first (resp. second) \textbf{tier numbers} (cf. \cite{Saito1985}). 
\vskip 0.2in

\noindent{\fbox{\textbf{An isolated case}}} \hskip 0.1in 
The ERS of type $BB_2^{\vee (2)\ast}$ is defined as follows:
\[ R=(R(BC_2)_s+\Z a+\Z b) \cup (R(BC_2)_m+L_{0,0}) \cup (R(BC_2)_l+ 2 \Z a+2 \Z b). \]

\subsubsection{$\diamond$-type}\label{sect:intro_non-red-diamond}
The ERS of type $C^\vee C_l^{(2)\diamond}$ is defined as follows:
\[ R=(R(C^\vee C_l)_s+\Z a) \cup (R(C^\vee C_l)_m+\Z a) \cup (R(BC_l)_l+\{ \, ma+2nb\, \vert\,  m-n\equiv 0 [2]\, \}). \]

\subsection{Classification of mERSs with non-reduced affine quotient}
The classification of the reduced mERSs $(R,G)$ with non-reduced  $R/G$ is given as follows:
\begin{thm}\label{thm_classification-reduced} 
Let $(R,G)$ be a reduced marked elliptic root system with non-reduced affine quotient $R/G$ belonging to a real vector space
with a symmetric bilinear form of signature $(l,2,0)$. Then it is isomorphic to one of the following:
\[ BC_l^{(1,2)}, \qquad BC_l^{(1,1)\ast}, \qquad BC_l^{(4,2)}, \qquad BC_l^{(4,4)\ast},  \qquad  BC_l^{(2,2)\sigma}(2) \qquad l\geq 1\]
and 
\[ BC_l^{(2,2)\sigma}(1) \qquad l\geq 2. \] 
\end{thm}
It is clear that by switching $a$ and $b$, 
from the ERSs of type $BC_l^{(1,2)}$, $BC_l^{(4,2)}$, $BC_l^{(2,2)\sigma}(1)$ and $BC_l^{(2,2)\sigma}(2)$, we obtain the ERSs of type $BC_l^{(2,1)}$, $BC_l^{(2,4)}$, 
$BC_l^{(2,2)}(1)$ and $BC_l^{(2,2)}(2)$ respectively, discovered by K. Saito in \cite{Saito1985}. The ERSs of type $BC_l^{(1,1)\ast}$ and $BC_l^{(4,4)\ast}$ has been considered by S. Azam in \cite{Azam2002}. Hence, Saito's list together with these $2$ ERSs due to Azam provide us a complete list of reduced ERSs. \\

As for the non-reduced mERSs $(R,G)$, the classification is given as follows: 

\begin{thm}\label{thm_classification-non-reduced}
Let $(R,G)$ be a non-reduced marked elliptic root system belonging to a real vector space
with a symmetric bilinear form of signature $(l,2,0)$.  Then, it is isomorphic to one of the following:
\begin{enumerate}
\item $R/G$ of type $BCC_l$ $(l\geq 1)$:
\begin{align*}
& BCC_l^{(1)}, \; BCC_l^{(1)\ast_0},\; BCC_l^{(1)\ast_{0'}}, \\
& BCC_l^{(2)}(1)\,  (l>1), \; BCC_l^{(2)}(2), \; BCC_l^{(2)\ast_p} \, (p \in \{0,1\}), \\
& BCC_l^{(4)}.
\end{align*}
\item $R/G$ of type $C^\vee BC_l$ $(l\geq 1)$:
\begin{align*}
& C^\vee BC_l^{(1)}, \\
& C^\vee BC_l^{(2)}(1)\,  (l>1), \; C^\vee BC_l^{(2)}(2), \; C^\vee BC_l^{(2)\ast_p} \, (p \in \{0,1\}), 
\\
& C^\vee BC_l^{(4)}, \; C^\vee BC_l^{(4)\ast_0},\; C^\vee BC_l^{(4)\ast_{0'}}.
\end{align*}
\item $R/G$ of type $BB_l^\vee$ $(l\geq 2)$:
\begin{align*}
&BB_l^{\vee \, (1)}, \; BB_l^{\vee \, (2)}(1), \;  BB_l^{\vee \, (2)}(2), \; BB_l^{\vee \, (4)}, \\
&BB_2^{\vee \, (2)\ast}.
\end{align*}
\item $R/G$ of type $C^\vee C_l$ $(l\geq 1)$:
\begin{align*}
&C^\vee C_l^{(1)}, \; C^\vee C_l^{(1)\ast_p} \, (p \in \{0,1\}), \\
&C^\vee C_l^{(2)}(1)\,  (l>1), \; C^\vee C_l^{(2)}(2), \; C^\vee C_l^{(2)\ast_s}, \; C^\vee C_l^{(2)\ast_l}, \\
&C^\vee C_l^{(2)\ast_{0}}, \; C^\vee C_l^{(2)\ast_1}, \; 
C^\vee C_l^{(2)\ast_{1'}}, \; C^\vee C_l^{(2)\diamond}, \\
&C^\vee C_l^{(4)}, \; C^\vee C_l^{(4)\ast_p} \, (p \in \{0,1\}).
\end{align*}
\end{enumerate}
\end{thm}
Here, we briefly explain how to classify mERSs with non-reduced quotient. Let $(R,G)$ be such a mERS of rank $l\geq 2$. (The rank 1 case can be handled by direct computations.) Set $R_+=R_m \cup R_l$ and $R_-=R_m \cup R_s$. It can be shown that the mERSs $(R_\pm,G)$ are of $C_l$-type (resp. $B_l$-type). By the classification theorem 
(cf. \cite{Saito1985}) of mERSs with reduced quotient, we have a list of possible mERSs $(R_\pm, G)$. In particular, the only possible mERS $(R_m,G)$ is of type $D_l^{(1,1)}$, and $R(D_2)+L_{0,0}$ for $l=2$ in addition. Here, $D_2$ (resp. $D_3$) is viewed as $A_1\times A_1$ (resp. $A_3$). Realizing $R_\pm$ in $F$ so that $R_+ \cap R_-=R_m=R(D_l^{(1,1)})$, it suffices to \textit{glue} them to obtain $R$. Namely, rotating one of $R_\pm$, say $R_+$, by an automorphism $\varphi \in \Aut(R_m)\subset GL(F)$, it is enough to verify that $R_- \cup \varphi(R_+)$ becomes a generalized root system. It turns out that isomorphism classes of such root systems can be parametrized by the double coset
\[ \Aut(R_+) \backslash \Aut(R_m)\slash \Aut(R_-). \]
\medskip
In \cite{AKY}, the authors described $BC_l$-type root systems as
\[ R(S,L,E)=(R(BC_l)_s+S) \cup (R(BC_l)_m+L) \cup (R(BC_l)_l+E) \]
for some discrete subsets $S, L, E \subset \rad_\Z (I)$ and gave certain conditions on the triple $(S, L, E)$ for $l>1$ and the pair $(S,E)$ for $l=1$, following the idea developed in \cite{AABGP}. They have thus classified the triples and pairs via some combinatorial studies on them, whereas our argument is based on the structure of automorphism groups of ERSs of $BCD$-type.

\subsection{Classification of non-reduced ERSs}
Theorem \ref{thm_classification-non-reduced} provides the isomorphism classes of mERS with non-reduced affine quotient. But, as root systems, some mERSs are isomorphic. Indeed, 
\begin{thm} Among the above $35$ non-reduced marked elliptic root systems, we have the following isomorphisms as root systems:
\begin{enumerate}
\item Via the isomorphism $a \leftrightarrow b$:
\begin{align*}
&R(BCC_l^{(2)}(1)) \cong R(BB_l^{\vee (1)}), \qquad 
  R(C^\vee BC_l^{(2)}(1)) \cong R(BB_l^{\vee (4)}), \\
&R(BCC_l^{(2)}(2)) \cong R(C^\vee C_l^{(1)}), \qquad 
  R(C^\vee BC_l^{(2)}(2)) \cong R(C^\vee C_l^{(4)}), \\
&R(BCC_l^{(2)\ast_p}) \cong R(C^\vee C_l^{(1)\ast_p}) \qquad
  R(C^\vee BC_l^{(2)\ast_p}) \cong R(C^\vee C_l^{(4)\ast_p}) \qquad (p \in \{0,1\}), \\
&\phantom{ABCDEFGHI} R(BB_l^{\vee (2)}(2)) \cong R(C^\vee C_l^{(2)}(1)).
\end{align*}
\item Via the isomorphism $a \mapsto a+b, \; b \mapsto b$:
\begin{align*}
&R(BCC_l^{(1)\ast_{0'}}) \cong R(BCC_l^{(1)\ast_0}), \qquad
  R(C^\vee BC_l^{(4)\ast_{0'}}) \cong R(C^\vee BC_l^{(4)\ast_0}), \\
&\phantom{ABCDEFGHI} R(C^\vee C_l^{(2)\ast_{1'}}) \cong R(C^\vee C_l^{(2)\ast_1}).
\end{align*}
\item Via an exotic isomorphism $\&$ $a \leftrightarrow b$:
\[ R(BCC_l^{(4)}) \cong R(C^\vee BC_l^{(1)}) \cong R(C^\vee C_l^{(2)\diamond}). \]
Indeed, the exotic one is given by
 \begin{align*}
 R(C^\vee C_l^{(2)\diamond})=
 &\left( R(BC_l)_s+\Z(a+2b)+ \Z b\right) 
 \cup \left( R(BC_l)_m+\Z(a+2b)+ 2\Z b\right) \\
  \cup &\left( R(BC_l)_l+\Z(a+2b)+4\Z b\right) 
 \cong R(C^\vee BC_l^{(1)}).
 \end{align*}
\end{enumerate}
\end{thm}

Thus, we only have $21$ isomorphic classes of non-reduced ERSs. 
\bibliographystyle{plain}
\bibliography{nr-mers}

\begin{thebibliography}{1}

\bibitem{AABGP}
B.~N. Allison, S.~Azam, S.~Berman, Y.~Gao, and A.~Pianzola.
\newblock Extended affine {L}ie algebras and their root systems.
\newblock {\em Mem. Amer. Math. Soc.}, 126(603):x+122, 1997.

\bibitem{Azam2002}
S.~Azam.
\newblock Extended affine root systems.
\newblock {\em J. Lie Theory}, 12(2):515--527, 2002.

\bibitem{AKY}
S.~Azam, V.~Khalili, and M.~Yousofzadeh.
\newblock Extended affine root systems of type {BC}.
\newblock {\em J. Lie Theory}, 15(1):145--181, 2005.

\bibitem{Kac1969}
V.~G. Kac.
\newblock Automorphisms of finite order of semisimple {L}ie algebras.
\newblock {\em Funkcional. Anal. i Prilo\v{z}en.}, 3(3):94--96, 1969.

\bibitem{Macdonald1972}
I.~G. Macdonald.
\newblock Affine root systems and {D}edekind's {$\eta $}-function.
\newblock {\em Invent. Math.}, 15:91--143, 1972.

\bibitem{Moody1969}
R.~V. Moody.
\newblock Euclidean {L}ie algebras.
\newblock {\em Canadian J. Math.}, 21:1432--1454, 1969.

\bibitem{Saito1974}
K.~Saito.
\newblock Einfach-elliptische {S}ingularit\"{a}ten.
\newblock {\em Invent. Math.}, 23:289--325, 1974.

\bibitem{Saito1985}
K.~Saito.
\newblock Extended affine root systems. {I}. {C}oxeter transformations.
\newblock {\em Publ. Res. Inst. Math. Sci.}, 21(1):75--179, 1985.

\end{thebibliography}
\newpage
\phantom{ABC}
\vskip 5in

{\small
\begin{tabular}{@{}l@{}}%
\textsc{A. Fialowski, }\\
\text{Department of Computer Algebra,} \\
\text{E\"{o}tv\"{o}s Lor\'{a}nd University,} \\
\text{P\'{a}zm\'{a}ny P\'{e}ter sét\'{a}ny 1$\slash$C,} \\
\text{H-1117 Budapest, Hungary} \\ 
\textit{E-mail address}: \texttt{fialowsk@inf.elte.hu, alice.fialowski@gmail.com}
\end{tabular}

\vskip 0.5in
\begin{tabular}{@{}l@{}}%
\textsc{K. Iohara, }\\
\text{Univ Lyon, Universit\'{e} Claude Bernard Lyon 1,} \\ 
\text{CNRS, Institut Camille Jordan UMR 5208} \\
\text{F-69622 Villeurbanne, France} \\ 
\textit{E-mail address}: \texttt{iohara@math.univ-lyon1.fr}
\end{tabular}

\vskip 0.5in

\begin{tabular}{@{}l@{}}%
\textsc{Y. Saito, }\\
\text{Department of Mathematics,} \\ 
\text{Rikkyo University,} \\
\text{Toshima-ku, Tokyo 171-8501, Japan} \\ 
\textit{E-mail address}: \texttt{yoshihisa@rikkyo.ac.jp, yosihisa@ms.u-tokyo.ac.jp}
\end{tabular}

}

\end{document}